\documentclass[10pt]{amsart}
\usepackage{amssymb, amsthm, amsmath}
\usepackage[all]{xy}
\usepackage{graphicx}
\usepackage{psfrag}

\newtheorem{conj}{Conjecture}

\theoremstyle{definition}

\newcommand\C{{\mathbb C}}

\newcommand\cH{{\mathcal H}}

\newcommand\Z{{\mathbb Z}}

\newcommand\al{\alpha}

\newcommand\la{\lambda}
\newcommand\s{{\sigma}}

\newcommand\eqto{\stackrel{\lower1.5pt\hbox{$\scriptstyle\sim\,$}}\to}

\newcommand{\skipline}{\vspace{\baselineskip}}

\DeclareMathOperator{\Sp}{Sp}

\DeclareMathOperator{\GL}{GL}

\DeclareMathOperator{\IG}{IG}

\DeclareMathOperator{\G}{G}

\DeclareMathOperator{\HH}{\mathrm{H}}

\begin{document}

\title[The roots of mathematical thought]
{The roots of mathematical thought}

\date{October 16, 2020}

\author{Harry~Tamvakis} \address{University of Maryland, Department of
Mathematics, William E. Kirwan Hall, 4176 Campus Drive, 
College Park, MD 20742, USA}
\email{harryt@math.umd.edu}

\maketitle

Six years ago, I wrote the opinion piece \cite{T2}, titled
`Mathematics is a quest for truth'. My aim in this sequel is to expand
upon the sentiment expressed there that `the roots of mathematical
thought lie within the deepest recesses of the human mind, where {\em
  logos} and {\em mythos} come together, in search of our very
nature'. I will examine the moment when new mathematics is created,
which is also a moment of discovery and revelation. There is a certain
mystical quality to this event, which is to a large extent very
personal. In conclusion, I will discuss what I believe mathematical
research can tell us about ourselves, and our role in the world.

I agree with Dieudonn\'e \cite[II.6]{D} and Hardy \cite{H} that the
main reason which compels us to do research in mathematics is
intellectual curiosity, the attraction of enigmas, the desire to know
the truth. Note however that neither \cite{T2} nor the present
paper claim that anything mathematical is actually true.

To give credence to my thoughts, I have to draw upon my own
experience, which inevitably entails some self-promotion. I apologize
in advance for the latter, and will have something more to say about
that at the end. I feel compelled to write because I have not seen
this topic discussed in the same way before, despite what I strongly
suspect, that I am not alone. Eleven years ago, I solved my favorite
problem in all of mathematics, in a manner so effortless that I did not
realize the full significance of the event until much later. The story
of how that happened will accompany the more important points to be
made along the way.

\skipline

I was reading books by my second birthday, and exhibited curiosity
about all sorts of things. Around that time, my family moved to the
United States from my native country of Greece, so I had two languages
to play with. An incident that affected me occurred in fourth grade
elementary school, in Skokie, Illinois, when our teacher taught us
about area, and showed how to derive the formula $\pi R^2$ for the
area of a circle. I remember the proof to this day. He cut the circle
into two semicircles, and divided each of them using radii into a
large and equal number of triangular wedges. The two halves of the
circle unfold along these lines by straightening their outer rim, like
a sliced orange, and fit together to form an (approximate) rectangle
whose side lengths are $\pi R$ and $R$, respectively.

I noticed that there was one shape our teacher had mentioned without
giving us a formula for its area: the `oval'. When I asked `What is
the area of an oval?', he said it requires calculus, and I would have
to wait until college, or perhaps late in high school, to learn
that. This led me on a quest to teach myself calculus and find the
answer, which lasted for years. The plot thickened when we moved back
to Greece a few years later, as I did not know the translation of
Latin terms like `calculus' and `oval' into Greek. The few books I
managed to gather on the subject were not very helpful, and I never
understood the precise definition of a limit until it was taught to me
in class, during my last year of high school.

In fact, a proof of the formula for the area of an ellipse does not
require calculus.  An ellipse with semiaxes of lengths $a$ and $b$
with $a<b$ can be defined as the intersection of a cylinder, whose
horizontal cross section is a circle $C$ of radius $a$, with a slanted
plane. The ellipse is obtained by dividing $C$ into two halves and
stretching each half in the vertical direction by a factor of $b/a$.
The area of the ellipse is thus $b/a$ times the area $\pi a^2$ of $C$,
or $\pi ab$.  This argument using similar triangles is not beyond the
comprehension of a curious fourth grader. However I am glad that
my teacher did not tell me this, because the journey I embarked on to
discover the answer was much more valuable than the destination.

\skipline

In secondary school, I found the answers that mathematics gave to the
key question of `why?' more satisfying than those offered in my
science classes. Although mathematics and science are close relatives
and cross-pollinate each other, there remains a fundamental distinction
between them, in that scientific validity depends upon agreement with
experiments, whereas mathematical facts rely on logically rigorous
proofs from abstract first principles. As such, the latter are among
the very few things that (so far) survive the test of time.

This proximity to truth and beauty was a major part of the attraction
of mathematics for me. Another was my participation and success in
math competitions, which exposed me to challenging problems. In
downtown Athens, there were weekly informal training sessions for
mathematical olympiads, organized by the Hellenic Mathematical Society,
and I decided to go have a look. The very first problem I faced upon
my arrival there was the following.

\medskip

{\em A school has 1000 students. Each student has a locker, and these
  are numbered from 1 to 1000. The first student opens all the
  lockers. The second student closes every second locker. The third
  student goes to every third locker, closing those lockers which are
  open, and opening those which are closed.  This continues in the
  same manner, with the $n$th student going to every locker whose
  number is a multiple of $n$, and opening/closing those which are
  closed/open, for each $n\leq 1000$. How many lockers will be open at 
  the end of this procedure?}

\medskip

After overcoming my initial reaction, here is how I solved this:
Consider a fixed locker, say locker number 12.  The students touching
this locker are numbers 1, 2, 3, 4, 6, and 12, which are the divisors
of 12. The locker remains closed at the end because 12 has an 
{\em even} number of divisors. Observe that the divisors of any positive
integer $n$ can be organized in pairs of the form $(d,n/d)$, such as
the pairs $(1,12)$, $(2,6)$, and $(3,4)$ in our example. Therefore $n$
will have an even number of divisors {\em unless} there is a
divisor $d$ of $n$ such that $d = n/d$, that is, unless $n=d^2$ is a
square. Hence the lockers numbered $1, 4, 9, 16, 25, \ldots $ will
remain open, and, since there are $31$ squares less than $1000$, the
answer to the problem is $31$.

This delightful conundrum, whose solution requires only basic
arithmetic, has two important lessons to teach us. The first is the
utility of working on a specific example, and then trying to
generalize from there. The second lies deeper. When confronted with
the question, one's first impression is that of a long line of
lockers, some of them changing their state with every
passing student. It is a picture of dizzying complexity,
similar to the famous sieve of Eratosthenes, which detects the
primes. In order to solve the problem, one needs to concentrate on a
fixed locker and carefully analyze what happens to it. This involves
changing your point of view, from a {\em global/horizontal} to a {\em
  local/vertical} one (see the figure below).

\medskip

\[
\includegraphics[scale=0.50]{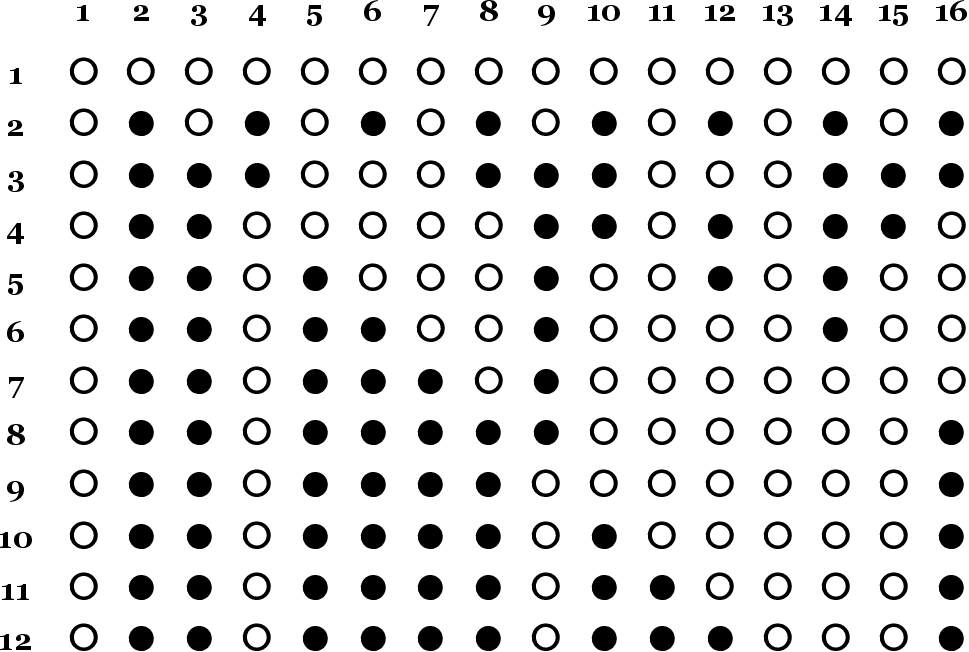}
\]

\medskip

The importance of a new perspective shedding light in an area of
mathematics will be a theme of this narrative. The examples are
myriad, at every level of the subject, for instance the introduction
of variable quantities and instantaneous rates of change, which turned
the ancient problem of the computation of areas from a static to a
dynamic one, and led to the fundamental theorem of calculus.

\skipline

My undergraduate study at the University of Athens provided me with a
solid background in analysis, which impressed my teachers during my
first year of graduate school at the University of Chicago.  It came
as a surprise to me that five years later, I was writing my thesis in
arithmetic algebraic geometry. My favorite subject was complex
analysis, and under the direction of my first advisor Narasimhan, I
went from the study of Riemann surfaces and Hodge theory to
Grothendieck's theory of schemes and intersection theory in algebraic
geometry. Narasimhan then suggested that I have a look at Gillet and
Soul\'e's arithmetic intersection theory, a generalization of Arakelov
theory to higher dimensions, which added number theory to the
mix. Since the glue that held all this together was intersection
theory, it made sense after my second topic exam to switch advisors to
Fulton, who was also at Chicago, and whose brilliant book on the
subject had kept me going.

One of the main difficulties with higher dimensional Arakelov theory
is a lack of examples where explicit computations are
possible. Maillot had succeeded in computing arithmetic intersection
numbers on Grassmannians, and I decided to try to extend his work to
more general flag manifolds. This involved moving beyond the case of
hermitian symmetric spaces, which have a canonical theory of harmonic
differential forms, and also learning a fair amount of combinatorics.

The advice I received from Fulton was to first try doing some small
examples. This is certainly the way one should begin any such project,
however I found to my frustration that I was unable to compute even
{\em one} new example using the known methods. Instead, what was
required was a different perspective: my Ph.D.\ dissertation was one
of the first works to use {\em representation theory} in Arakelov
geometry, and led the way to a new {\em arithmetic Schubert
  calculus}. During the long road to its completion, I came to
appreciate all that mathematics has to offer, and the beauty found
in areas that seem far from one's initial interests.

I remember the moment when I realized that I could solve my thesis
problem.  I was alone in my living room, sitting on the couch, when my
mind turned again, as it had so many times before, to the calculation
I sought to achieve. And then, in a sudden flash of inspiration, I saw
how I could do it. The solution combined in a novel way ingredients
that I had developed gradually over the previous years, and proved
that all the natural intersection numbers were rational numbers. On
that occasion, I experienced the ineffable, incomparable feeling of
exciting new mathematics taking shape right before (or rather, behind)
my eyes.

The following excerpt from a famous poem by Dionysios Solomos
describes in the best way that I know the exhilaration that the
researchers of the unknown experience at the moment of mathematical
discovery.\footnote{I am indebted to the Greek philosopher Christos
  Malevitsis, who loved this passage and discussed its significance
  in his writings.  The translation and interpretation here are my own.}

\medskip

\medskip

\qquad \qquad {\em Mother, magnanimous in suffering and glory,

\qquad \qquad Though your children forever live in hidden mystery

\qquad \qquad In reflection and in dream, what grace have these eyes,

\qquad \qquad These very eyes, to behold thee in the deserted forest} $\ldots$

\medskip

\medskip

The `Mother' is the ontological mother -- the mother that gives birth
to us all.  She is the {\em truth} -- no matter how defined -- that
mathematicians seek in their research endeavors. Her children live in
hidden mystery, alone and in the dark, searching for a pathway to that
truth. We have two essential tools at our disposal: deductive
reasoning (reflection) and guesses or conjectures, which are rooted in
human intuition and fantasy (dreams). Much of the work happens
subconsciously, over long periods of time.  And when the moment of
enlightenment finally comes, one is left speechless in admiration and
awe.  The only feeling then is one of wonder, joy, and gratitude for
the gift that was bestowed upon you, to be present there with your
mind's eye open, at the very moment when the Mother reveals a bit more
of her light. We have here an essentially otherworldly experience, as
the researcher understands that he or she is more the recipient or
channel of knowledge, rather than its creator. No other area of
inquiry illustrates in such a profound way the illusion of human
agency: we apparently have original mathematical ideas, but the
conclusions we obtain using them could not have been otherwise.

\skipline

There is ample evidence that new mathematics is discovered, often
independently by different people at different times, and is a
revelation to us. A striking example of this occurred in my first
joint paper with Kresch on the arithmetic Grassmannian. For such
spaces, Gillet and Soul\'e formulated arithmetic analogues of
Grothendieck's difficult standard conjectures on algebraic cycles.

Kresch and I used arithmetic Schubert calculus to prove these
conjectures for the Grassmannian of lines in projective space
\cite{KT}. In the course of our work, we were led to make a conjecture
of our own about a certain family of {\em Racah polynomials}. This
indirect and miraculous connection between the problem we set out to
solve, the theory of hypergeometric orthogonal polynomials, and the
Racah coefficients or $6$-$j$ symbols in quantum mechanics came as a
complete surprise.

To state our conjecture, let $k$, $\ell$, and $n$ be integers, and define
\begin{equation}
\label{Req}
R(k,\ell,n) := \sum_{i=0}^{n-1}(-1)^i
{k \choose i} {k+i \choose i} {\ell \choose i} {\ell+i \choose i}
{n-1 \choose i}^{-1} {n+i \choose i}^{-1}.
\end{equation}

\begin{conj}
\label{conj2}
For any integers $k,\ell,n$ with $0\leq k,\ell\leq n-1$, we have 
\[
-1\leq  R(k,\ell,n)\leq 1.
\]
\end{conj}

\medskip
We found much computer evidence in support of Conjecture \ref{conj2},
and proved it when $\ell\leq 3$ or $\ell=n-1$. The latter case is
interesting: although one can show, using the Wilf-Zeilberger method,
that in general there is no `closed form' for the sum defining
$R(k,\ell,n)$, when $\ell=n-1$ there is such a formula, namely
\begin{equation}
\label{PfSa}
R(k,n-1,n) = \sum_{i=0}^k(-1)^i\frac{n}{n+i}\binom{k}{i}\binom{k+i}{i}
=\frac{(1-n)(2-n)\cdots (k-n)}{(1+n)(2+n)\cdots (k+n)}.
\end{equation}
The second equality in (\ref{PfSa}) is a special case of an identity
proved by Pfaff, 220 years ago, and (independently!) by Saalsch\"utz,
130 years ago.

The most vivid evidence for Conjecture \ref{conj2} is the following
picture, which was kindly provided to us by Wilf. In the figure, we
let $n:=51$ and plot the values $R(k,\ell,51)$ on the lattice points
$0\leq k,\ell\leq 50$, then connect the resulting dots by line
segments. In fact, the bound $|R(k,\ell,51)| \leq 1$ fails badly if $k$
and $\ell$ are taken to be {\em real} parameters lying in $[0,50]$. In
other words, the conjecture depends essentially on the {\em integrality} 
of $k$ and $\ell$.

\medskip

\[
\includegraphics[scale=0.31]{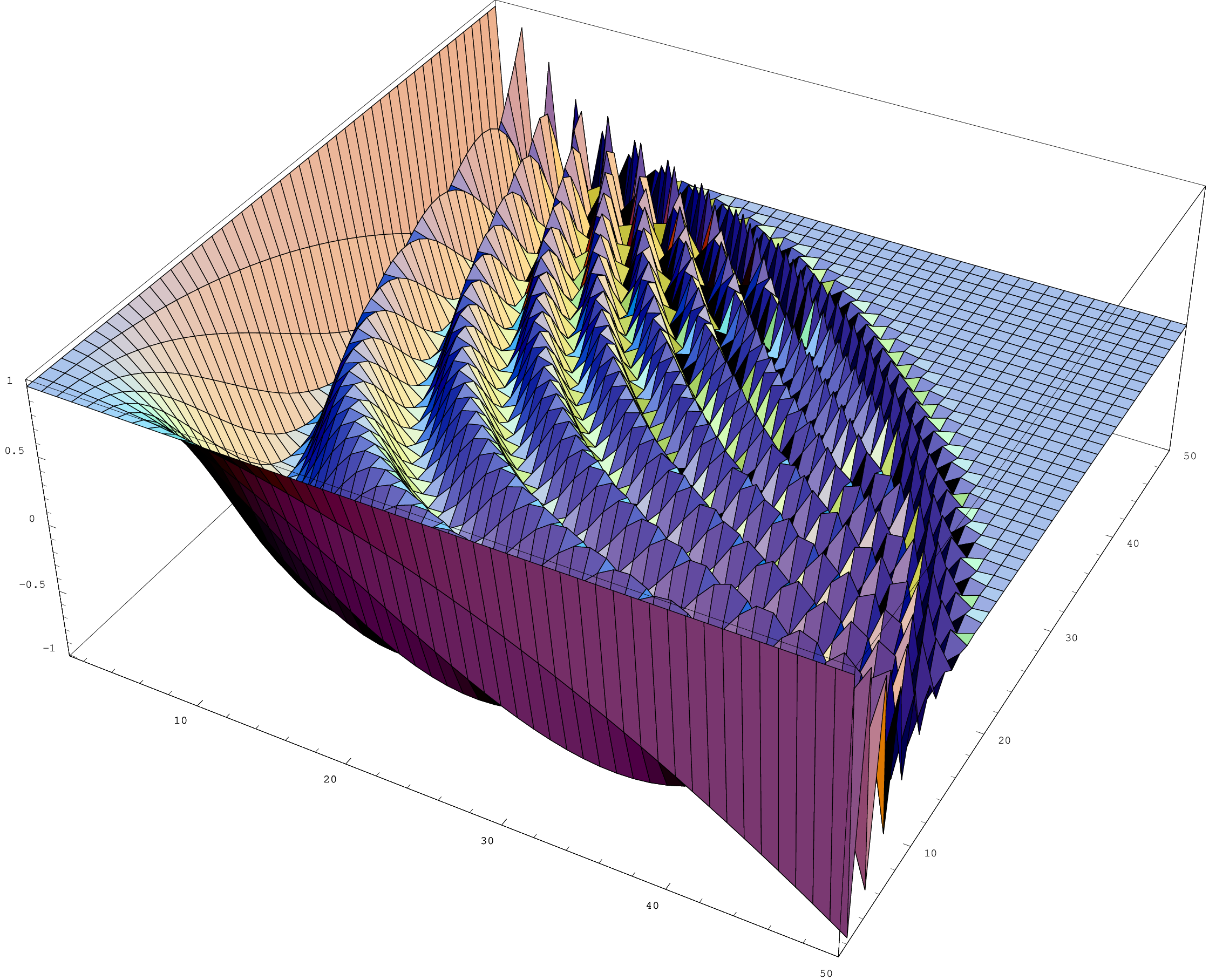}
\]

\medskip
Conjecture \ref{conj2}, which surprised Askey, remains open as
of this writing.  Fortunately for Kresch and me, the arithmetic Hodge
index conjecture for the Grassmannian $\G(2,n)$ is equivalent to the
weaker inequality
\[
\sum_{k=1}^{n-1} (-1)^{k+1} R(k,\ell,n) \cH_k <
\sum_{k=1}^{n-1} \cH_k
\]
which was within our reach. Here
$\cH_k:=1+\frac{1}{2}+\cdots+\frac{1}{k}$ is a {\em harmonic number}.

\skipline

My research on the Arakelov theory of Grassmannians and flag manifolds
was followed by further joint papers on their quantum cohomology
rings. These exotic intersection theories are both deformations of the
classical Schubert calculus. The latter is the study of the usual
cohomology ring of the same spaces, and an old and rich subject.  The
work of twentieth century mathematicians provided algorithms to do
computations in Schubert calculus, but I realized that there remained
fundamental gaps in our understanding.  A perplexing mystery beyond
the Lie type A was the problem of {\em representing polynomials}, to
which we turn next.

We begin with the example of the Grassmannian $\G(m,n)$, consisting
of all $m$-dimensional $\C$-linear subspaces of $\C^n$. This space has
a natural decomposition 
\begin{equation}
\label{Ceq}
\G(m,n)=\bigcup_I X_I
\end{equation}
into {\em Schubert cells} $X_I$, one for each subset $I:=\{i_1,
\cdots,i_m\}$ of $\{1,\ldots,n\}$ with $|I|=m$. Every subspace $V$ in
$\G(m,n)$ can be represented uniquely by an $m\times n$ matrix $A$ in
reduced row echelon form with row space $V$. We have $V\in X_I$ if and
only if the pivot `1's in $A$ lie in columns $i_1,\ldots,i_m$. The
{\em Schubert class} $[X_I]$ is the cohomology class of the closure of
$X_I$. The cell decomposition (\ref{Ceq}) implies that the Schubert
classes form an additive basis for the cohomology group of $\G(m,n)$:
\[
\HH^*(\G(m,n),\Z) = \bigoplus_I\Z[X_I].
\]

To understand the ring structure of $\HH^*(\G(m,n),\Z)$, it is better
to parametrize the Schubert classes by integer partitions
$\la=(\la_1\geq\cdots\geq \la_m)$ with $\la_1\leq n-m$ rather than by
subsets $I$. The formula $\la_r:=n-m+r-i_r$ for $1\leq r \leq m$ gives
a bijection between these two parameter spaces. If $I$ corresponds to
$\lambda$, we denote the Schubert class $[X_I]$ by $\s_\la$.  For
every integer $p\in [1,n-m]$, the class $c_p:=\sigma_{(p,0,\ldots,0)}$
is known as a {\em special Schubert class}.  Giambelli \cite{G} was
able to express a general Schubert class $\s_\la$ as a polynomial in
the special classes $c_1,\ldots,c_{n-m}$:
\begin{equation}
\label{Geq}
\s_\la = \det(c_{\la_i+j-i})_{1\leq i,j\leq m}.
\end{equation}
In the Giambelli formula (\ref{Geq}), and in later equations, we set
$c_0:=1$ and $c_p:=0$ if $p\notin [0,n-m]$.  The polynomials on the
right hand side of (\ref{Geq}) are algebraic representatives for the
Schubert classes, and lead naturally to an algebraic model for the
cohomology ring of Grassmannians.

In fact, the determinant in equation (\ref{Geq}) had appeared in the
work of Jacobi and Trudi, 60 years earlier. The symmetric functions
studied by these authors were eventually named {\em Schur
  polynomials}, in honor of Schur's work relating them to the
characters of the general linear group. Although this remarkable
connection between the cohomology of Grassmannians and representation
theory was (and continues to be) influential, later research showed
that it breaks down when one looks at more general homogeneous spaces.

In the twentieth century, the study of the Grassmannian and related
symmetric spaces was enriched by writing them as quotients of Lie
groups. The general linear group $\GL_n=\GL_n(\C)$ acts transitively
on $\G(m,n)$, and the stabilizer of the subspace $\langle e_1,\ldots,
e_m\rangle$ under this action is the maximal parabolic subgroup $P$ of
invertible matrices in the block form
\[
\left(
\begin{array}{c|c}
* & * \\  \hline
0 & *
\end{array}
\right)
\]
where the lower left block is an $(n-m)\times m$ zero matrix. It
follows that $\G(m,n)=\GL_n/P$.  The Weyl group of $\GL_n$ is the
symmetric group $S_n$, and if $B\subset P$ denotes the Borel subgroup
of upper triangular matrices, then we have a decomposition
\[
\GL_n/P = \bigcup_{w \in S_n^P}Bw P/P,
\]
where $S_n^P$ is the set of permutations $w$ such that $w(i)<w(i+1)$
for all $i\neq m$. This agrees with the Schubert cell decomposition
(\ref{Ceq}) of $\G(m,n)$ given earlier.

The above picture generalizes to the case when $G$ is a classical
complex Lie group and $P$ is any parabolic subgroup of $G$, so that
the quotient space $G/P$ is a compact manifold parametrizing
(isotropic) partial flags of subspaces in $\C^N$. If $B\subset P$ is a Borel
subgroup, then we have a cell decomposition
\[
G/P= \bigcup_{w \in W^P}Bw P/P
\]
and a corresponding direct sum decomposition
\begin{equation}
\label{Heq}
\HH^*(G/P,\Z) = \bigoplus_{w\in W^P} \Z \,\sigma_w
\end{equation}
of the cohomology group of $G/P$, where $W^P$ is a certain subset of
the Weyl group of $G$.  The {\em Giambelli problem} is to determine
the analogue of formula (\ref{Geq}), that is, to find canonical
polynomial representatives for the Schubert classes $\s_w$ on $G/P$.

Once again, the answer required a change in perspective. The most
challenging and original part of this was joint work with Buch and
Kresch, which examined the case of symplectic Grassmannians. Here we
equip $\C^{2n}$ with a non-degenerate skew-symmetric bilinear form
$(\ ,\,)$, and say that a linear subspace $V$ of $\C^{2n}$ is {\em
  isotropic} if the restriction of $(\ ,\,)$ to $V$ vanishes
identically. We let $\IG=\IG(n-k,2n)=\Sp_{2n}/P$ denote the symplectic
Grassmannian consisting of all isotropic subspaces of dimension $n-k$,
for some fixed $k\geq 0$. The Schubert classes $\s_\la$ on $\IG$ can
be indexed by {\em $k$-strict partitions} $\la$. The condition
`$k$-strict' means that all parts $\la_i$ of $\la$ with $\la_i>k$ are
distinct, and reflects the fact that the rows of the matrices in
reduced row echelon form which represent $V$ must be pairwise
orthogonal.

Using a Pieri rule for the products $c_p\cdot \s_\la$ of a special
Schubert class $c_p:=\sigma_{(p,0,\ldots,0)}$ with a general one, we
proved that the $c_p$ for $p\in [1,n+k]$ generate the ring
$\HH^*(\IG,\Z)$, and could access the Schubert calculus
there. With the help of a computer, we observed that known
determinantal and Pfaffian formulas represented $\s_\la$ in extreme
cases, when all the parts of $\la$ were at most $k$, or,
respectively, greater than $k$. However, a search of the extensive
literature for an operation on matrices which interpolates naturally
between a determinant and a Pfaffian in the sense required proved
fruitless. To make matters worse, there was no a priori reason why
there should be {\em any} nice formula for $\s_\la$, and, given
the presence of relations among the $c_p$, there were plenty of
reasons to be pessimistic.  A situation all too familiar to those
researching the unknown: we were groping in the dark.

Instead of the old language of determinants and Pfaffians, the answer
we found to the Giambelli problem for $\IG$ employed Young's {\em
  raising operators}. Given an integer sequence
$\al:=(\al_1,\al_2,\ldots)$ with finite support and indices $i<j$,
define
\[
R_{ij}(\alpha) := (\alpha_1,\ldots,\alpha_i+1,\ldots,\alpha_j-1,
\ldots).
\] 
A raising operator $R$ is any monomial in the $R_{ij}$'s.  Moreover,
let $c_\al:=c_{\al_1}c_{\al_2}\cdots$, and for any raising operator
$R$, let $R\,c_{\al} := c_{R(\al)}$.  The main theorem of \cite{BKT}
states that for any $k$-strict partition $\la$, we have
\begin{equation}
\label{Theq}
\s_\la = \prod_{i<j} (1-R_{ij})\prod_{\la_i+\la_j > 2k+j-i}
(1+R_{ij})^{-1} \, c_\la
\end{equation}
in $\HH^*(\IG,\Z)$, 
where the first product is over all pairs $i<j$ and second product is
over pairs $i<j$ such that $\la_i+\la_j > 2k+j-i$. 

If the partition $\la$ 
satisfies $\la_i\leq k$ for all $i$, then (\ref{Theq}) becomes
\begin{equation}
\label{Veq}
\s_\la = \prod_{i<j} (1-R_{ij})\, c_\la = \det(c_{\la_i+j-i})_{i,j}
\end{equation}
so we recover the Giambelli formula (\ref{Geq}) in this case. The second equality in 
(\ref{Veq}) is a formal consequence of the Vandermonde identity.
For example, we have
\[
(1-R_{12})\, c_{(a,b)} = c_{(a,b)} - c_{(a+1,b-1)} = 
c_ac_b-c_{a+1}c_{b-1} = \left|\begin{array}{cc}
c_a & c_{a+1} \\ c_{b-1} & c_b
\end{array}\right|.
\]
At the other extreme, if $\la_i>k$ for all non-zero parts $\la_i$ of $\la$, then 
(\ref{Theq}) becomes
\[
\s_\la = \prod_{i<j}\frac{1-R_{ij}}{1+R_{ij}}\, c_\la = 
\mathrm{Pfaffian}\left(\frac{1-R_{12}}{1+R_{12}}\, c_{\la_i,\la_j} \right)_{i<j}.
\]
Here the second equality follows from a classical result due to Schur.

Since Young introduced them in the 1930s, the raising operators
$R_{ij}$ made occasional appearances, notably in the theory of
Hall-Littlewood functions, and in the book \cite{M}. However, for both
historical and mathematical reasons, they were hardly ever used. Our
paper \cite{BKT} was the first to show that raising operators play an
essential role in geometry, in the Giambelli formulas for isotropic
Grassmannians, and to employ them in their proofs. It was a lengthy
process, with many ups and downs, from the moment in 2003 when formula
(\ref{Theq}) was conjectured, until our proof of it was complete. The
raising operator expressions in (\ref{Theq}), which depend on the
indexing partition $\lambda$, do not enjoy the same alternating
properties as determinants, and we had to learn how to work with them
from scratch.

We called the polynomial on the right hand side of (\ref{Theq}) a {\em
  theta polynomial}, and its even orthogonal counterpart an {\em eta
  polynomial}.  The theta and eta polynomials are the analogues of the
Schur polynomials in the symplectic and orthogonal Lie types, for the
purposes of Schubert calculus, and of geometry more generally. The
expository paper \cite{T3} gives further information about this
correspondence.

\skipline

The work \cite{BKT}, along with several companion papers, was
announced on arXiv in the fall of 2008. The following summer, I sensed
that all the required tools were in place to settle the Giambelli
problem in the general case (\ref{Heq}).  I wrote the paper \cite{T1}
very quickly, during the first few weeks of July 2009. The solution
combined many different strands of prior research, and also solved the
problem when $G/P$ varies in an algebraic family, to obtain formulas
for the cohomology classes of {\em degeneracy loci of vector
  bundles}. The main theorems gave unique, combinatorially explicit,
and {\em intrinsic} Chern class formulas for all the degeneracy loci
involved.

I finished writing \cite{T1} and stared in disbelief at the end
result. Although it did require some new ingredients, the paper seemed
almost trivial, with most proofs consisting of just a few lines. The
article \cite{T1} introduced a new, intrinsic point of view in
Schubert calculus, one that came entirely naturally to me, the
culmination of an understanding of the subject that had developed over
many years.

We are not always so fortunate that the answers to our hardest
questions turn out to be so simple. Nevertheless, I believe that
simplicity is the hallmark of truth. I also began to appreciate,
quoting E.\ Artin, that our difficulty is not in the proofs, but in
learning {\em what} to prove. Indeed, the ability to ask the right
questions is critical, and something that can only be taught by
experience.

To the casual observer looking at the sequence of papers leading up to
\cite{T1}, it appears as though the author had it all planned out,
assembling the necessary components over time, until the final
synthesis. Of course, this is completely false, as a detailed
examination of the record shows. Many of the pieces were found while
solving problems in Arakelov theory and quantum cohomology, not
directly related to \cite{T1}. Moreover, a large part of the work was
intuitive, with some sections included in papers not because they were
needed there, but because they were too beautiful to omit, and trusting
there would be an application someday.

\skipline

As I began my education, I was drawn to mathematics because I love
truth, and sought it there. I soon realized that although mathematics
is seemingly humanity's most credible attempt at finding permanence
and truth, it has not reached that goal, not by a long shot. Nor does
mathematics have answers yet to our deepest and most pressing
existential questions.  Since aquaintance with truth is the basis of
all knowledge, I conclude with Socrates that in essence, {\em I know
  nothing}. Everything written below should therefore be taken with
that caveat in mind.

Because I can't honestly claim to possess knowledge with any
certainty, what I am left with are beliefs. Like any scientist, my
beliefs are supported by evidence culled from my life experiences.
The aforementioned ones are the most relevant for the purposes of this
brief exposition, but many more remain unsaid.

I believe that the apparent existence and consistency of the mental
structure we call mathematics is a miracle, as wondrous as the
universe around us. The beautiful equations such as (\ref{Req}) and
(\ref{Theq}) which emerge out of our collective mind are as
unforeseeable as any revelatory event. These miracles are unlike those
that depend on religious faith: once seen, they can be reproduced,
shared, and admired together with other human beings. However, they
are far too surprising and otherworldly to be explained away as our
inventions, or a product of solely rational thought.

I am under no illusion that I knew what to expect, or even that I was
in the driver's seat, when \cite{T1} was written. Like most of my
research efforts, the work began with a pen, paper, and wishful
thinking. Although the article soon materialized on my desk, I have no
satisfactory explanation as to why all the different ingredients
required were there just when they were needed, and fit together like
a charm, so that the proofs turned out to be so easy.

To paraphrase Gibran, it is not up to you to direct the course of
mathematics, for {\em mathematics, if it finds you worthy, directs
  your course}.  I therefore cannot in good faith take credit for
\cite{T1}, or by extension any of my papers. This is not just because
of the debt that \cite{T1} owes, like all scientific research, to the
many other works by various authors that preceded it. More
importantly, I am convinced that the fact that \cite{T1} was going to
be written was {\em determined in advance}.

This is what my career has helped to teach me about the human
condition: that the choices we make, consciously or not, are
predetermined. Science supports the notion that our sense of self and
freedom of the will are both illusory. Mathematicians are in a unique
position to understand this, because their success depends on
circumstances which are figments of their imagination, and yet
assuredly beyond their control. Moreover, in both the research process
and its fruits, mathematicians experience and perceive miracles. At
the same time, as a poster in my daughter's elementary school points
out, math is everywhere.

Assuming that human society can reconcile with determinism and
overcome its anthropocentric posture without destroying itself in the
process, what is left for the human race to strive for? Recall the
importance of adopting the right perspective, and also the Mother,
magnanimous in suffering and glory. In my opinion, what remains is to
collectively move closer to Her point of view. This will require a
humble reckoning with the truth about our endeavors, an admission of
ignorance, and above all, bestowing this attitude towards life to our
children, and looking to them -- for forgiveness first, and,
ultimately, for guidance.

\end{document}